\documentclass{article}

\usepackage[utf8]{inputenc} 
\usepackage[T1]{fontenc}    
\usepackage{url}            
\usepackage{booktabs} 
\usepackage[top=2cm, bottom=2cm, left=4cm, right=4cm]{geometry}
\usepackage{float}
\usepackage{stmaryrd}
\usepackage{amsfonts}
\usepackage{amssymb}
\usepackage{amsthm}
\usepackage{amsmath}
\usepackage{mathabx}
\usepackage{bigints}
\usepackage{algpseudocode}
\usepackage{bbm}
\usepackage{graphicx}
\usepackage{subcaption}
\usepackage{booktabs}
\usepackage{multirow}
\usepackage{longtable}
\usepackage[export]{adjustbox}
\usepackage{ragged2e}
\usepackage{threeparttable}     
\usepackage{nicefrac}       
\usepackage{microtype}      
\usepackage{xcolor}         
\usepackage{amsmath}       
\def\bx{\mathbf{x}}

\def\VV{\mathbb{V}}

\def\bx{\mathbf{x}}
\def\bX{\mathbf{X}}

\newcommand{\UU}{\mathbb{U}}

\newcommand{\R}{I\!\!R}
\newcommand{\N}{I\!\!N}

\newcommand{\Pp}{\mathbbm{P}}
\newcommand{\E}{\mbox{E}}

\newcommand{\1}{\mathbbm{1}}

\newtheorem{theorem}{Theorem}
\newtheorem{definition}{Definition}

\newtheorem{corollary}{Corollary}

\newtheorem{example}{\noindent Example}

\usepackage[linesnumbered,ruled,vlined]{algorithm2e}
\SetKwInput{KwRes}{Result}%
\SetKwIF{Si}{SinonSi}{else if}{if}{then}{else if}{else}{end}%
\SetKwFor{Tq}{While}{do}{end}
\SetKwFor{PourCh}{For each}{do}{end}

\title{Computing conservative  probabilities of rare events with surrogates}
\author{%
  Nicolas Bousquet\thanks{SINCLAIR AI Laboratory, EDF R{\&}D  \& LPSM, Sorbonne Université} 
}
\date{}

\begin{document}

\maketitle

\begin{abstract}
This article provides a critical review of the main methods used to produce conservative estimators of probabilities of rare events, or critical failures, for reliability and certification studies in the broadest sense. These probabilities must theoretically be calculated from simulations of (certified) numerical models, but which typically suffer from prohibitive computational costs. This occurs frequently, for instance, for complex and critical industrial systems. We focus therefore in adapting the common use of surrogates to replace these numerical models, the aim being to offer a high level of confidence in the results. We suggest avenues of research to improve the guarantees currently reachable. \\

{\bf Keywords}: rare events, probability of failure, surrogates, accelerated Monte Carlo, probability upper bound, reliability, safety, guarantee, kriging, polynomial chaos expansion, neural network
\end{abstract}


\section{Introduction}

The estimation of probabilities of a rare event (often described as a critical failure) defined by 
\begin{eqnarray*}
p & := & P(g(X)< y) 
\end{eqnarray*}
where  $g:\Omega=[0,1]^d \to \R$ is a deterministic function, $y\in\R$ such that $p>0$ and $X$ follows a uniform distribution on $\Omega$, is the basic principle of many structural reliability studies \cite{lemaire2013structural}. Actually, $X$ can be given more complicated distributions, but this framework remains very general.  Difficulties occur when $p$ can be very low and $g$ can only be known through knowledge of a small number of realisations of $Y=g(X)$, either by simulation of $X$, or from an observed sample of realizations of $X$ collected in a possibly sequential manner. Simulation-based approaches are given particular consideration in what follows, as they provide a theoretical underpinning to purely empirical approaches. In this framework, many methods have been proposed to produce  statistical estimators of $p$ based on $n$ queries  $x\mapsto g(x)$, with reduced asymptotic variance compared to the standard Monte Carlo (MC) approach, providing the well-known estimator  $p_n$ such that
\begin{eqnarray*}
\sqrt{n}\left(p_n - p \right) & \xrightarrow[{\cal{L}}]{n\to\infty} {\cal{N}}\left(0,\sigma^2\right).
\end{eqnarray*}
where $\sigma^2=p(1-p)$
See \cite{morio2014survey} for a recent survey of such techniques. Their best representatives are  based on  importance sampling (IS) \cite{bucklew2004introduction}, possibly fed by tools from large deviation theory in high-dimensional settings \cite{tong2023large,schorlepp2023scalable}, sequential IS \cite{PAPAIOANNOU201666}  and other adaptive IS methods as the cross-entropy method \cite{de2005tutorial}, directional sampling \cite{morio2015estimation}, line sampling \cite{lu2008reliability,de2015advanced,papaioannou2021combination}, multilevel Monte Carlo \cite{elfverson2016multilevel,hajiali2021adaptive}
, variational approaches \cite{valsson2014variational,friedli2023energybased} 
or importance splitting \cite{LEC07} (e.g., subset simulation \cite{au2001estimation,au2003subset} and adaptive multilevel splitting \cite{cerou2012sequential,brehier2015analysis,cerou2019adaptive}). If the best asymptotic variance that one may expect from these latter approaches is $p^2\log(1/p)\ll \sigma^2$ \cite{guyader2011simulation}, they require sampling multiple Markov chains, which can dramatically increase the number of queries $x\mapsto g(x)$ \cite{bernard2021recursive}. For this reason, many other approaches are taking the gamble of building statistical surrogates (or meta-models) $x\mapsto \hat{g}_m(x)$  from the characteristics of $x\mapsto g(x)$ and a design of experiments (DOE) ${\bf x_m}=x_1,\ldots,x_m\in\Omega^m$, to reduce the overall cost of the calculation. Methods described in next paragraphs are summarized in Table \ref{tab:methods}.  

\subsection{Surrogate-based methods}

A first class of methods considers a nonintrusive statistical approximation of $x\mapsto g(x)$ through quadratic response surfaces \cite{GAYTON200399} or, much more commonly, in a Bayesian framework, through almost surely continuous Gaussian processes (kriging-based regression) whose mean defines the deterministic surrogate $\hat{g}_m$ \cite{santner2003design}. These imply assumptions of regularity (especially on the correlation structure \cite{rasmussen2010gaussian}). 
Intrusive (adjoint-based) techniques can lead to define surrogates as \emph{reduced order models (ROM)}\footnote{As recalled by \cite{cerou2023adaptive}, such methods 
are known to be "very efficient for the numerical approximation of problems involving the repeated solution of parametric Partial Differential Equations (PDEs)".}, ie. defined using a reduced order basis \cite{quarteroni2015reduced}, that may come with error bounds with respect to $x\mapsto g(x)$, as considered in \cite{cerou2023adaptive}. 

Coming back to kriging-based surrogates, various sequential simulation techniques can then be used to benefit from the meta-model uncertainty, for example based on targeted Integrated Mean Square Error (tIMSE) \cite{picheny2010adaptive} or Bayesian Stepwise Uncertainty Reduction (SUR) strategies \cite{bect2012sequential}, possibly combined with subset sampling methods \cite{bect2017bayesian}. Such strategies are defended by strong theoretical results \cite{bect2019supermartingale}. Over the years, many adaptive variants of kriging-based methods have been proposed to estimate $p$ (e.g. \cite{bichon2008efficient}). Active learning strategies, mixing sampling within $\Omega$ (MC, other better space filling techniques as Latin Hypercube Sampling, importance sampling,  subset sampling, etc.) and kriging to produce estimates of $p$, have been particularly popularized. Let us cite the so-called AK-MCS methods (see \cite{lelievre2018ak} and \cite{moustapha2022active} for recent reviews, and references therein). In engineering, kriging-based approaches appear useful to determine the subset of $\Omega$ leading to failures (e.g., \cite{azzimonti2021adaptive,marrel2022icscream}), a dual gain in estimating $p$. Most famous competitors to Gaussian process meta-modeling are polynomial chaos expansions  \cite{schobi2015polynomial}, implied within similar combinations of techniques to reach the estimation of very low probabilities. Both approaches have numerous merits in the general field of uncertainty quantification (e.g., versatility, easiness for conducting sensitivity analyses) and are the subject of major research aimed at reducing their computational cost in storage and inference (typ. $\mathrm{O}(m^2)$ and $\mathrm{O}(m^3)$ for kriging) and handling large dimensions $d$. Finally, the nesting of transformations in the regression function, for example through deep Gaussian processes \cite{damianou2013deep}, makes it possible to capture the complexities of the topological manifold defined by $x\mapsto g(x)$, at the cost of a large sample size. Neural networks and their universal approximation capability provide a final subclass of surrogates used in the latter context, aiming at diminishing the regularity assumptions. See \cite{lelievre2018ak} and \cite{jakeman2022surrogate} for two recent reviews of the various surrogate techniques used in the field of reliability, to complete this brief summary.  \\  

A second class of methods is concerned solely with the construction of a surrogate $\hat{\Gamma}_m$ of the failure (or \emph{limit state}) surface 
\begin{eqnarray*}
\Gamma & = & \left\{x\in\Omega, \ g(x)=y \right\}
\end{eqnarray*}
which is defined under continuity assumptions, and considering that the problem is a perfectly separable binary classification problem. This generalizes in finding excursion sets, as described hereinafter. 
Engineering methods like FORM/SORM \cite{ditlevsen1996structural} make the assumption that $\Gamma$ can be locally approximated by linear or quadratic surfaces, and transform the estimation problem into an optimization problem. They are usually combined with IS, as techniques based on large deviation theory \cite{tong2023large}. Combined with subset simulation, Support Vector Machines (SVM) or combinations of SVM were proposed by \cite{bourinet2011assessing} and \cite{bousquet2018approximation} to approximate $\Gamma$ (see \cite{roy2023support} for a review), while neural networks were preferred by \cite{papadrakakis2002reliability,bousquet2012accelerated,lieu2022adaptive}. Kriging-based classification was proposed by \cite{dubourg2011reliability}, while \cite{LI20108966} chose polynomial chaos expansion, accompanied with cross-entropy importance sampling \cite{li2011efficient}.  These approximations are often combined, again, with sequential sampling strategies. The AK-MCS methods evoked hereinbefore fall into this category too, as they are underlyingly based on a classifier (approximating how a vector $x$ is far from $\Gamma$)   defined by a kriging predictor.   

\subsection{Ensuring conservatism}

When it comes to the real-life estimation of $p$ to characterize the safety of critical systems (e.g., nuclear/aeronautical/spatial  structures and processes), where $g$ is typically a certified numerical model, 
the question arises of the credit to be given to estimators based on  $\hat{g}_m$ or $\hat{\Gamma}_m$. Do they provide "good", "sufficient" estimators of the true probability $p$? Can we establish insurance rules based on these estimators?
Despite the continuous improvement of these methods and their combinations, and the subsequent reduction of the error attributable to the use of surrogates\footnote{A part of this error being linked to dimension-reduction, see \cite{jakeman2022surrogate}.}  -- in particular by adaptive sampling techniques and selected queries $x\mapsto g(x)$ \cite{au2007augmenting,grooteman2008adaptive,li2011efficient} -- most surrogate-based methods are still lacking of theoretical guarantees adapted to the context. The asymptotic convergence to $p$ of any estimator of the kind
\begin{eqnarray}
\hat{p}_{n,m} & = & \frac{1}{n} \sum\limits_{i=1}^n \1_{\{\hat{g}_m(x_i) < y\}} \ \ \ \ \ \ \text{\it (MC or IS with surrogate of $g$)} \label{MC.with.surrogate.g}
\end{eqnarray}
or
\begin{eqnarray}
\hat{p}_{n,m} & = &  \frac{1}{n} \sum\limits_{i=1}^n \1_{\{x_i \prec \hat{\Gamma}_m\}}, \ \ \ \ \ \  \text{\it (MC or IS with surrogate of $\Gamma$)} \label{MC.with.surrogate.Gamma}
\end{eqnarray}
where $\prec$ is some partial ordering rule (allowing binary classification), 
appears \emph{weak}, in the sense it does not answer to the requirements that can be expected from certification authorities\footnote{As said in \cite{Ducoffe2020}, ``Certification authorities will most likely require safeguards", while the authors of \cite{jakeman2022surrogate} insist: ``The expectation is inappropriate for risk-averse stakeholders, as it does not capture any notion of variability and does not quantify rare events"}. While the consistency of (\ref{MC.with.surrogate.g}) and (\ref{MC.with.surrogate.Gamma}) in $(m,n)$ is of course necessary, this asymptotic property remains a figment of the imagination when $g$ is expensive, and/or $d$ is large.  And, more generally, because we can hardly define the reality of an asymptotic regime, even with a central limit theorem.  Stronger guarantees should be  non-asymptotic, for instance provided through concentration inequalities (or PAC-Bayes-type inequalities) as
\begin{eqnarray}
P\left(\left|\hat{p}_{n,m} - p \right|>\varepsilon\right) & \leq & \alpha_{m,n,\varepsilon} \label{concentration.inequality}
\end{eqnarray}
where $\alpha_{m,n,\varepsilon}$ remains ideally very low for achievable $(m,n)$ (at the order of $p$) and small $\varepsilon$ and $\alpha_{m,n,\varepsilon}\to 0$ when $(m,n)\to\infty$ for any $\varepsilon$. 

Other strong guarantees are likely to come from the specific choice of ROM, as they often come with error controls with respect to $g$. Evoked in the recent work \cite{cerou2023adaptive}, such properties seem not to have been yet too much explored, maybe because  such ROM are typically made to represent regular behaviors rather than extreme ones. 
Finally, other strong guarantees are related to robustness to variations on $g$ (or, similarly, on the real distribution of $X$, through a transport from the uniform distribution introduced hereinbefore). Even most of MC-based methods provide estimates of $p$ that lack of robustness, as their coefficient of variation goes to infinity when $p$ goes to $0$ \cite{glynn2009robustness}. This type of problem involves both robustness studies and uncertainty quantification. Robust analysis theories struggling against the probabilistic misspecification of epistemic uncertainties, such as the info-gap theory \cite{ben2006info,ben2019info}, which uses convex epistemic models\footnote{But can be associated to other epistemic models (e.g. possibility distributions, p-boxes, Dempster-Schafer structures...).}, are used in {\it hybrid} reliability studies \cite{wang2010reliability}. 
\begin{eqnarray*}
X \ \sim \ \underbrace{\text{distribution $P$}}_{uncertain} & \xrightarrow[]{Consequence} & p\in(p^-_n,p^+_n)
\end{eqnarray*} 
Using tools derived from random set theory \cite{tonon2004using}, they can be used to define bounds on probabilities of failure. See  \cite{ajenjo2022info} for more information.  However, computing these bounds is itself subject to the same computational problems as that of a single probability $p$ for which the model $X\sim \Pp$  is well defined, and usually requires reduction variance techniques or/and surrogates. An alternative is to make geometric assumptions on the distribution of the output $g(X)$ to get conservative bounds, in the spirit of \cite{strigini2014bounds}. These authors provide rules in the specific case of survival analysis, where prior information on lifetime distributions is available. \\

One way of simultaneously addressing the   previous concerns is to produce {\it conservative estimators} of $p$, ie. upper bounds $p^+_{n,m}\geq p$, such that $p^+_{n,m}\to p$ and $p^+_{n,m}-p$ be not too wide \cite{bousquet2012accelerated}.   Ideally, one would like the order of magnitude of $p^+_{n,m}$ to be similar to that of $p$: if $p\sim 10^{-q}$, then it is hoped to get $10^{-q}\leq p^+_{n,m}\ll 10^{-q+1}$. Note that if some algorithm can be run to determine $p^+_{n,m}\in[p,1)$, by symmetry the same algorithm can be  adapted to determine a lower bound $p^-_{n,m}\in(0,p]$, provided $p>0$. \\

Noticing that $p=\mbox{Vol}(\Omega_{y})$ where $\Omega_{y}$ is the so-called \emph{excursion set} (or \emph{level set}) \cite{bolin2015excursion,azzimonti2021adaptive}
\begin{eqnarray*}
\Omega_{y} & = & \left\{x\in \Omega: \ g(x)<y\right\},
\end{eqnarray*}
any upper bound $p^+_{n,m}$ for $p$ corresponds to the volume of a subset  $\Omega^+_{y,m,n}$ such that $\Omega_y\subseteq \Omega^+_{y,m,n}\subsetneq\Omega$. The complementary subset $\widebar{\Omega}^+_{y,m,n}=\Omega/\Omega^+_{y,m,n}$ corresponds to a set smaller than the \emph{safe} set $\widebar{\Omega}_y=\{x\in\Omega: \ g(x)\geq y\}$. When the elements of $\widebar{\Omega}^+_{y,m,n}$ are included in $\widebar{\Omega}_y$ with a large probability $\beta$, such sets $\widebar{\Omega}^+_{y,m,n}$ are called \emph{conservative} by \cite{french2013spatio,bolin2015excursion}. Still in a Bayesian setting, these latter authors proposed a kriging-based approach to estimate this conservative set (thanks to the posterior distribution of $g$) from a fixed DOE. This work was recently extended by \cite{azzimonti2021adaptive}, who proposed adaptive strategies based on specific SUR criteria to reduce the uncertainty of this estimator. While offering excellent performance in terms of volume recovery on low-dimensional cases, this approach costs $\sim$ 200 queries $x\mapsto g(x)$ on a real industrial example in 2D. While it has not been used specifically to compute upper bounds, it partially answers to the questions above, as it is conditioned by the relevance of the choice of the Gaussian process.  \\

To the best of our knowledge, obtaining appropriate guarantees for the certification of a calculation of a probability of failure using surrogates therefore remains a tough problem. The aim of this article is to summarize some useful results for further research in this area. First (Section \ref{sec:geometry}), we recall some strong conservatism results linked to the use of certain geometric properties of $g$. These are difficult to check by nature, but can be much more easily used to select certain classes of surrogates. They provide theoretically deterministic bounds on a probability $p$ calculated from a surrogate. 
We also consider in Section \ref{sec:trials}  approaches inspired by \cite{Ducoffe2020} and \cite{jakeman2022surrogate}, which seek to "bias" a surrogate to ensure a  form of conservatism in the estimation of $p$. A discussion of future prospects concludes this work, focusing on the obstacles to be overcome.

\begin{table}[hbtp]
\centering
\caption{A summary of most popular methods dedicated to the estimation of rare event probabilities.}
\label{tab:methods}

\vspace{0.2cm}

\begin{tabular}{lll}
\hline
\multicolumn{2}{l}{{\bf Variance reduction methods}} & { (survey: \cite{morio2014survey})} \\
\hline
Importance sampling (IS) & & \cite{bucklew2004introduction} \\
& + large deviation theory &\cite{tong2023large,schorlepp2023scalable} \\
& + spectral decomposition & \cite{yuan2021efficient} \\ 
& + variational approaches &\cite{valsson2014variational,friedli2023energybased} \\
& Sequential IS & \cite{PAPAIOANNOU201666} \\
& Cross-entropy IS & \cite{de2005tutorial} \\
Directional sampling & & \cite{morio2015estimation}  \\
Line sampling & & \cite{lu2008reliability,de2015advanced,papaioannou2021combination} \\
Multilevel Monte Carlo & & \cite{elfverson2016multilevel,hajiali2021adaptive} \\
Importance splitting & & \cite{LEC07} \\
& Subset simulation & \cite{au2001estimation,au2003subset} \\
& Adaptive multilevel splitting & \cite{guyader2011simulation,cerou2012sequential,brehier2015analysis,cerou2019adaptive} \\
& \\
\hline
\multicolumn{2}{l}{{\bf Methods based on a surrogate of $g$}} & {(survey: \cite{lelievre2018ak})} \\
\hline
& \\
Reduced order models & & \cite{quarteroni2015reduced} \\
Quadratic response surfaces & & \cite{GAYTON200399} \\
Mean predictors of  {Gaussian processes} & & \cite{santner2003design} \\
& + Sequential strategies  & \cite{bect2019supermartingale} \\
& Ex 1: Bichon criterion & \cite{bichon2008efficient} \\
& Ex 1: iMSE criterion & \cite{picheny2010adaptive} \\
& Ex 2: Stepwise Uncertainty Reduction (SUR) & \cite{bect2012sequential} \\
& Ex 3 : SUR + Subset sampling & 
\cite{bect2017bayesian} \\
& Ex 4 : Active learning (AK-MCS) & \cite{lelievre2018ak,moustapha2022active} \\
& Ex 5 : derived from inversion methods & \cite{marrel2022icscream} \\
& Ex 6 : derived from excursion sets & \cite{bolin2015excursion,azzimonti2021adaptive} \\
Polynomial chaos expansions & & \cite{schobi2015polynomial} \\
& \\
\hline
\multicolumn{2}{l}{{\bf Methods based on a surrogate of $\Gamma$}} & { (survey: \cite{roy2023support})} \\
\hline
& \\
Linear response surfaces (FORM) & & \cite{ditlevsen1996structural} \\
Quadratic response surfaces (SORM) & & \cite{ditlevsen1996structural} \\
& + IS & \cite{breitung2021sorm} \\
Support Vector Machines (SVM) & & \cite{bourinet2011assessing} \\
& Combination of SVM & \cite{bousquet2018approximation} \\
Neural networks & & \cite{papadrakakis2002reliability,bousquet2012accelerated} \\
& Adaptive strategies & \cite{lieu2022adaptive} \\
Adaptive kriging & & \cite{dubourg2011reliability} \\
Polynomial chaos expansion & & \cite{LI20108966} \\
& + cross-entropy IS & \cite{li2011efficient} \\
& \\
\hline
\end{tabular}
\end{table}
\section{Obtaining deterministic bounds on $p$ from geometrical properties}\label{sec:geometry}

The strongest desirable conservatism condition consists in producing deterministic upper bounds $p^+_{n}$ on $p$, independently on any additional $m-$sample that could be used to estimate a surrogate. The only known bounds are intrinsically linked to geometrical conditions, sometimes induced by conditions of regularity  which we recall and discuss below. For the sake of generality we denote by $\Pp$ the probability measure associated to $X$ over a subset of $\Omega$ (which, in our context, is simply the uniform measure). 

\subsection{Lipschitz smoothness}

\begin{figure}[hbtp]
\centering
\includegraphics[scale=0.4]{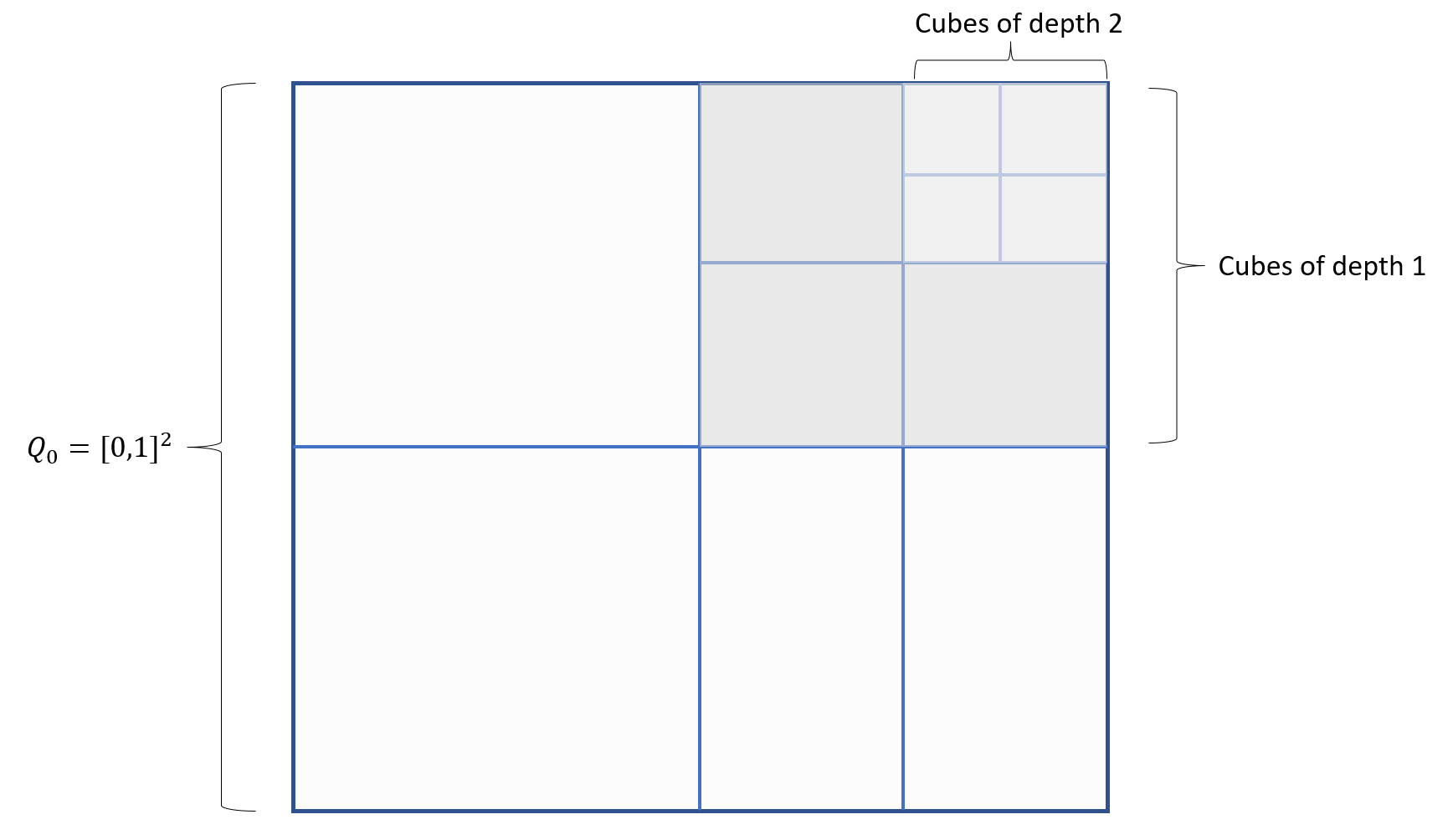}
\caption{An illustration of  various levels of partitioning $[0,1]^2$ by a family of dyadic cubes.}
\label{fig:dyadic.cubes}
\end{figure}

Theorem \ref{theo:bernard} is derived from a sequence of results recently obtained by 
\cite{bernard2021recursive}. It uses classical tools in harmonic analysis: dyadic cubes, for which a constructive definition must first be provided, illustrated by Figure \ref{fig:dyadic.cubes}.

\begin{definition}[Dyadic cubes.]
Dyadic cubes in $\Omega=[0,1]^d$ are a (possibly infinite) collection of cubes 
$$\left\{
\begin{array}{c}
Q_{0}, \\
Q_{1,1},\ldots,Q_{1,2^d}, \\
Q_{2,1,1},\ldots,Q_{2,1,2^d},Q_{2,1,1},\ldots,Q_{2,2,2^d},\ldots,\ldots,Q_{2,2^d,1},\ldots,Q_{2,2^d,2^d} \\
\ldots 
\end{array}
\right\}
$$
such that:
\begin{itemize} 
\item $Q_{j,\ldots}$ has sidelength $2^{-j}$ for $j\in\N$, with $Q_0=[0,1]^d$;
\item for any $j\in\N$, the $Q_{j,\ldots}$ define a partition of $\Omega$;
\item each $Q_{j,\ldots}$ has $2^d$ children cubes $Q_{j+1,\ldots}$ build by performing a $2^d-$ split of $Q_{j,\ldots}$; 
\item each $Q_{j,\ldots}$ has exactly one parent cube, for $j>0$.
\end{itemize}
A cube $Q$ with sidelength $2^{-j}$ is said to have a depth $j$.  
\end{definition}

Each cube $Q$ is then determined by its depth $j$ and its center, noted $c_Q$.  Reusing the notations of \cite{bernard2021recursive}, it can be labelled by a query to $x\mapsto g(x)$ as follows: 
\begin{itemize}
\item $Q\in {\cal{I}}$ if $g(c_Q)>y+L2^{-j-1}$ {\it (inside the safe subset of $\Omega$)},
\item $Q\in {\cal{O}}$ if $g(c_Q)<y-L2^{-j-1}$ {\it (outside the safe subset of $\Omega$)},
\item $Q\in {\cal{U}}$ otherwise.
\end{itemize}
Given a choice of maximal depth $k$, a recursive algorithm given in \cite{bernard2021recursive} allows $n$ cubes to be labelled, as precised in the following theorem. Denote then  ${\cal{I}}_n$ the set of all cubes labelled as ${\cal{I}}$ along this algorithm (and similarly let us denote ${\cal{O}}_n$ and ${\cal{U}}_n$). \\

\begin{theorem}[Derived from \cite{bernard2021recursive}.] \label{theo:bernard}
Under the following assumptions:
\begin{itemize}
\item[(i)] The distribution of $X$ on $\Omega$ admits a bounded density function with respect to the Lebesgue measure $\lambda$;
\item[(ii)] The function $g$ is assumed to be $L-$Lipschitz with respect to the supremum norm on $R^d$, ie.
\begin{eqnarray}
\left|g(x)-g(\tilde{x})\right| & \leq & L\|x-\tilde{x}\|_{\infty}, \ \ \ \ x,\tilde{x}\in\Omega;
\end{eqnarray}
\item[(iii)] There exists a constant $C>0$ such that
\begin{eqnarray}
\lambda\left(\left\{x\in\Omega: \ \left|g(x)-y\right|\leq \delta\right\}\right) & \leq & \frac{C}{L}\delta, \ \ \ \delta>0 \ \ \text{\it {(level set condition)}};
\end{eqnarray}
 \end{itemize}
then
\begin{eqnarray}
p^-_{n} \ \leq & p & \leq \ p^+_{n}
\end{eqnarray}
where
\begin{eqnarray*}
p^+_{n} \ = \ 1 - \sum\limits_{Q\in {\cal{I}}_n} \Pp(X\in Q) & \text{and} &
p^-_{n} \ = \ p^+_{n}-\sum\limits_{Q\in {\cal{U}}_n} \Pp(X\in Q),
\end{eqnarray*}
with, for $d\geq 2$, 
\begin{eqnarray*}
n & \leq & 4C 2^{d-1} \\
\text{and} \ \ \ p^+_{n} - p^-_{n} & \leq & C n^{-\frac{1}{d-1}},
\end{eqnarray*}
this latter result being the optimal convergence rate, in the sense it cannot be improved by any other algorithm defined under the sole general assumptions. 
\end{theorem}

In our context, Assumption (i) is always true, and note that $\Pp(X\in Q)=2^{-dj}$ for any $Q$ of sidelength $j$. Notice that \cite{bernard2021recursive} consider more broadly that $X$ can follow any distribution on $\Omega$ provided (i) remains true. In this case, the computation of the $\Pp(X\in Q)$ is not trivial and requires numerical methods, that bring stochastic errors on the estimation of bounds $(p^-_n,p^+_n)$. Using a splitting approach, the authors can still control the bounding of $p$ and the convergence of the difference of estimated bounds towards 0. To our knowledge, these results were only applied  to a univariate toy case with failure probability of the order $10^{-3}$ (see Table \ref{summary.proba}).

As explained by the authors, the level set Assumption (iii) reflects the fact that $g$ is not too much flat in the vicinity of the limit state (or level set, or failure surface) $\Gamma=g^{-1}(y)$. This condition appears mild under a continuous differentiability assumption on $g$. Besides, the knowledge of $M$ is not necessary for the algorithmic task. 

Estimating (or bounding) the Lipschtiz constant $L$ is however required and is a more serious difficulty for carrying out the real computation of the bounds. This is generally a difficult problem, as it involves being able to calculate a maximum distance between the outputs of $g$. However, if $g$ is a surrogate chosen as a neural network, owning the Lipschtiz property,  the tight estimation of $L$ has become over the recent years a challenge for which constructive solutions have recently been published. Especially, the approach proposed by \cite{fazlyab2019efficient}, that estimates $L$ by solving a semidefinite problem, has become very popular in the machine learning community. See besides \cite{pabbaraju2021estimating}  for the estimation of Lipschitz constants of monotone deep equilibrium models.

\subsection{Monotonicity}

A particular case of functions $x\mapsto g(x)$ are the monotonic ones, when monotonicity is understood with respect to the following partial order on $\R^d$ defining Pareto dominance :

\begin{definition}
\label{definition:PartialOrder}
Let $x=(x_1, \ldots,x_d)$ and $x' = (x'_1, \ldots, x'_d) \in \R^d$. 
If $x_i \leq x'_i$ for all $i= 1,\ldots,d$,   $x$ is said to be {\it dominated} by $x'$, and it is denoted $x \preceq_P x'$. 
\end{definition}

\begin{figure}[hbtp]
    \centering
    \includegraphics[scale=0.5]{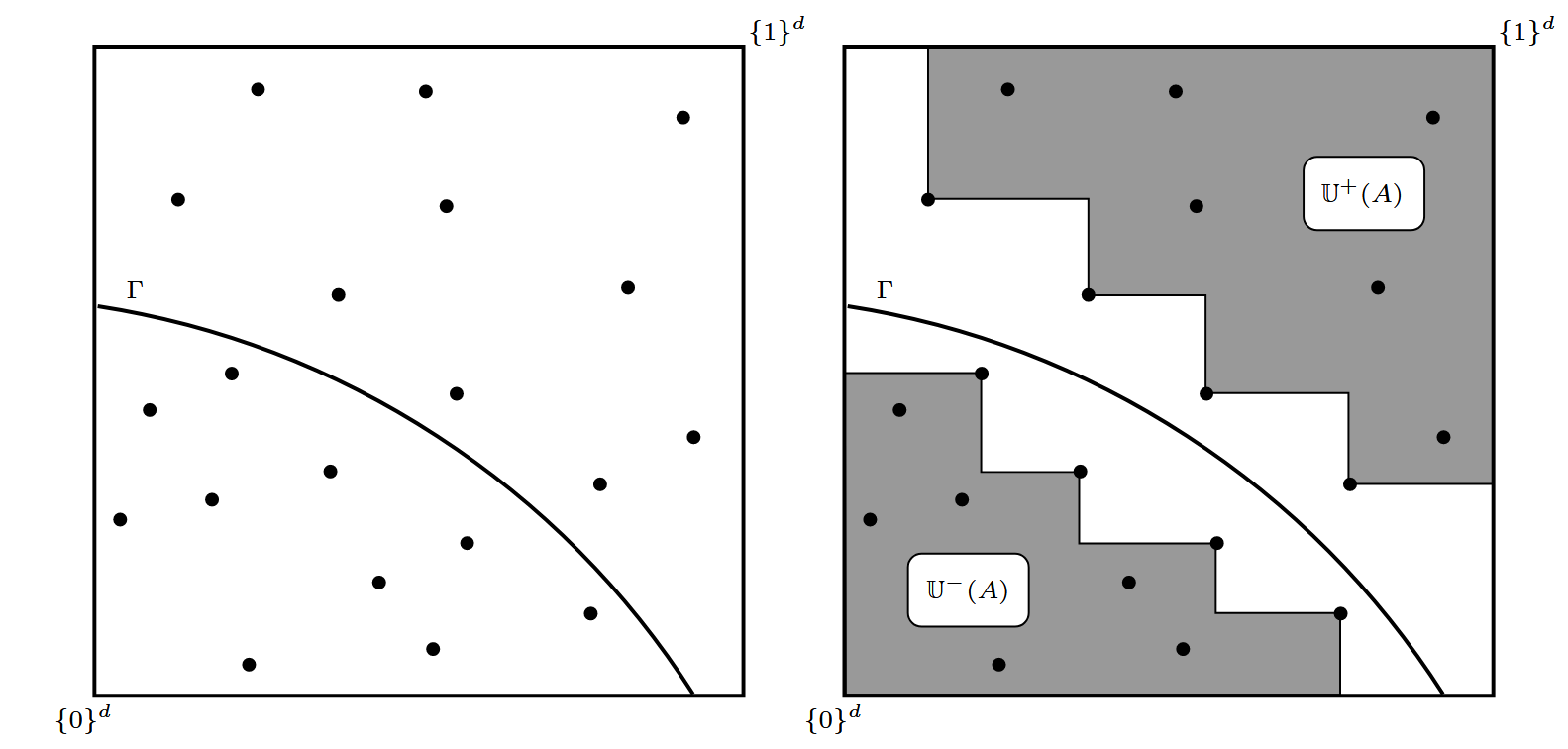}
    \caption{An illustration in dimension 2 of the bounding principle due to monotonicity. On the left is a set $A$ of $m$ points included in $[0,1]^d$. On the right, we plot the dominated spaces $\UU^-(A)$ and $\UU^+(A)$ surrounding $\Gamma$. It is clear that $\UU^-(A)\preceq_P \Gamma \preceq_P \UU^+(A)$. The volume of these subspaces allows to compute deterministic bounds for $p$: $\lambda(\UU^-(A))\leq p \leq 1 - \lambda(\UU^+(A))$, with $\lambda$ the Lebesgue measure.}
    \label{fig:mrm}
\end{figure}

To simplify things in the following, the monotonicity of $g$ is assumed to be similar for all dimensions in $\Omega$:  $g$ is supposed to be {\it globally increasing} (possibly at the price of a reparameterization), ie;  $\forall\{x,x'\} \in [0,1]^d\times [0,1]^d$ such that $x \preceq_P y'$, then $g(x) \leq g(x')$. An immediate consequence is that each $x\preceq_P x'$, when $g(x')\leq y$, is such that $g(x)\leq y$. This led \cite{de2009structural} and \cite{bousquet2012accelerated} then \cite{moutoussamy2013comparing} and \cite{bousquet2018approximation} to define a class of algorithms that simultaneously produce deterministic bounds $(p^-_{m},p^+_m)$ and consistent statistical estimators $\hat{p}_{n,m}$ of $p$, from the Lebesgue measures (volumes) of sequences of nested random dominated sets
\begin{eqnarray*}
\UU^-(A)  & = & \bigcup_{x \in A \cap \UU^-}\{u \in [0,1]^d: \ u \preceq x \}, \label{IS:eq:UUA:1}\\
\UU^+(A) & = &  \bigcup_{x \in A \cap \UU^+}\{u \in [0,1]^d: \ u \succeq \bx \}.\label{IS:eq:UUA:2} 
\end{eqnarray*}
with
\begin{eqnarray*}
\UU^- & = & \{x \in [0,1]^d: \ g(x) \leq y \} \ \text{and} \ 
\UU^+ \ = \ \{x \in [0,1]^d: \ g(x) > y \}.
\end{eqnarray*}
Figure \ref{fig:mrm} offers an 2D-illustration of the fundamentals of this method. The bounds obtained this way were used to guarantee the result of safety analysis for some industrial case-studies (e.g., \cite{sueur2013bounding}). Sequential surrogates of the limit state surface $\Gamma$ can be used to guide new design points to lower the volume $p^+_m-p^-_m$ (and jointly improve $\hat{p}_{n,m}$), but the core of the approach stands on the gain yielded by the deterministic (upper) bounds, which can be computed exactly in low dimensions \cite{bousquet2012accelerated} or else by specific numerical methods \cite{bousquet2018approximation}.  

Actually, for building such bounds it is enough to  that the limit state surface $\Gamma$ be $p-$monotone in the following sense:

\begin{definition}
\label{definition:AlphaMonotonic}
Denote $\mu$ the Lebesgue measure on $\R^d$. Let $A \subset [0,1]^d$. Define the dominated sets
\begin{eqnarray*}
\VV^-(A) & := & \bigcup_{x \in A}\{u \in [0,1]^d: \ u \preceq_P x \}, \label{Quantile:def:UUA:1}\\
\VV^+(A) & :=  & \bigcup_{x \in A}\{u \in [0,1]^d: \ u \succeq_P x \}. \label{Quantile:def:UUA:2}
\end{eqnarray*}
Let $\alpha \in (0,1)$ and $S$ be a set in $[0,1]^d$. 
The set $S$ is said $\alpha$-monotonic if for all $u,v \in S$ such that $u\neq v$, $u$ is not strictly dominated by $v$ and if $\mu(\VV^-(S)) = \alpha$.
\end{definition}

As a consequence of the monotonicity of $g$ (Proposition 5.1 in \cite{bousquet2018approximation}), a major feature of $\Gamma$  is its $p$-monotonicity, provided that $\Gamma$ is simply connected and $\lambda(\Gamma) = 0$. Note that these two last  assumptions are required in \cite{bousquet2018approximation} to estimate consistently $\Gamma$, interpreted as the frontier of a separable classification problem.

Monotonic sub-behaviors are often key to understand complex behaviors summed up by some simple rule as testing if $g(X)<y$, as demonstrated by Bshouty's monotone theory for learning Boolean functions \cite{bshouty1995exact},  many works related to monotonicity detection in such functions (e.g. \cite{743493}) and the importance of monotonicity in causal analysis and SAT theory.  Modeling and sensitivity analysis theories make extensive use of tools to detect monotonicity between inputs and $X$ and outputs $Y$ (typically via correlation coefficients). The exhibition of partial or total monotonicity properties of a model or phenomenon $x\mapsto g(x)$ defines a semantics allowing its interpretability and the obtaining of guarantees on the predictive outcomes $g(x)$ \cite{sharma2020higher}. Using monotonicity property from expert knowledge can help to produce more suitable surrogates in engineering prblems \cite{da2012gaussian,HAO2021101342,da2020gaussian}. Hence monotonicity  facilitates  the use of $g$ in critical domains such as clinical testing (e.g., \cite{mueller2023monotonicity}) or credit scoring (e.g., \cite{provenzano2020machine}). See the previous references for various industrial illustrations. 

For these reasons, numerous models that can be used as surrogates of $g$ or $\Gamma$ present global or local monotonicity properties, from the simplest one (linear model) to complex ones (e.g., deep lattice networks \cite{NIPS2017_464d828b} or GAMI-nets \cite{yang2021gami}, that allow {\it monotonization}). See for instance \cite{CANO2019168} for a review of monotonic classification. 

For instance, useful generalizations of linear models are the $m-$Continuous Piece-Wise Linear functions ($m-$CPWL)
which can be interpreted as  piecewise monotonic surrogates:

\begin{definition} \cite{Chua1988,ovchinnikov2000max,chen2022improved} A function $g:\Omega\to {\cal{Y}}$ is a $m-$Continuous Piece-Wise Linear function ($m-$CPWL) is there exists $K$ finite sets of disjoint complex polytopes ${A_k}_{k=1}^m$ such that $\bigcup_{k=1}^m A_k=\Omega$ and $g$ restricted to the domain $A_k$, denoted as $g_{|A_k}: A_k \ni x \mapsto g(x)$ is affine for each $k\in\{1,\ldots,m\}$. 
\end{definition}

This versatile class encompasses neural networks with piece-wise linear activations such as ReLu or hard tanh, that correspond to $\max(0,x)$ and $\max(-1,\min(1,x))$, respectively. As explained in \cite{Amoukou2023} (Chapter 3), {\it feedforward neural networks can be described as piece-wise linear functions that divide the input space into multiple linear , where the network itself behaves as an affine function within each region} \cite{hanin2019deep,hanin2019complexity,chen2022improved}. Reusing the notations in Theorem 2.2 from \cite{Amoukou2023}, we describe the hyperrectangles $A_k$ by 
\begin{eqnarray*}
A_k & = & \bigotimes\limits_{i=1}^p A_{i,k}
\end{eqnarray*}
where $A_{i,k}=[l_{i,k},r_{i,k}]$ with $(l_{i,k},r_{i,k})\in \widebar{\R}^2$. Each component $g_{|A_k}$ is represented by 
\begin{eqnarray*}
g_{|A_k}(x) & = & \sum\limits_{i=1}^p a_{i,k} x_i + b_k
\end{eqnarray*}
where the $(a_{i,k},b_k)$ are real numbers. Consequently, 
\begin{eqnarray}
g(x) & = & \sum\limits_{k=1}^m\left(\sum\limits_{i=1}^p a_{i,k} x_i + b_k\right) \1_{A_k}(x) 
\end{eqnarray}
which is clearly piecewise monotonic.


Numerical experiments conducted with the following toy example proposed in \cite{bousquet2012accelerated} were conducted until dimension 15 in \cite{bergere2021}. Approximate nested uniform sampling within the non-dominated space (defined by $\Omega/(\UU^+(A)\cup\UU^-(A))$ in Figure \ref{fig:mrm}) was produced using a semi-adaptive MCMC summarized in Appendix \ref{appendix:MCMC}.  

\begin{example}\label{toy.example.1}
Given a fixed dimension $d$, denote 
\begin{eqnarray*}
    V_d = \tilde{g}_d(Z)= Z_1/\big(X_1 + \sum_{i=2}^d Z_i \big)
\end{eqnarray*}
where the $Z_i$ are Gamma distributed: $Z_i\mathcal{G}(i+1,1)$, for $i=1,\ldots,d$.  The function $\tilde{g}_d$ is increasing in $Z_1$ and decreasing in all the variables $Z_i$ for $i \geq 2$. The same monotonicity can be deduced for $ x\mapsto \tilde{g}_d \circ T^{-1}(x)$ where $x\in\Omega$ and $T$ is obtained from the cumulative distribution functions of $Z_i$: $T=(F_{Z_1},\ldots,F_{Z_d})$.  The variable $V_d$ follows the Beta distribution $\mathcal{B}(2, 2^{-1}(d+1)(d+2)-3)$.
For $p \in [0,1]$, let $q_{d, p}^1 $ be the quantile of order $p$ of $Z_d$, and let $g_d(x) = \tilde{g}_d \circ T^{-1}(x) - q_{d, p} $. Hence, with $y=0$, $P(g_d(X)<0)=p$ for all $d\geq 2$.
\end{example}

\begin{figure}[hbt!]
\centering
\includegraphics[scale=0.6]{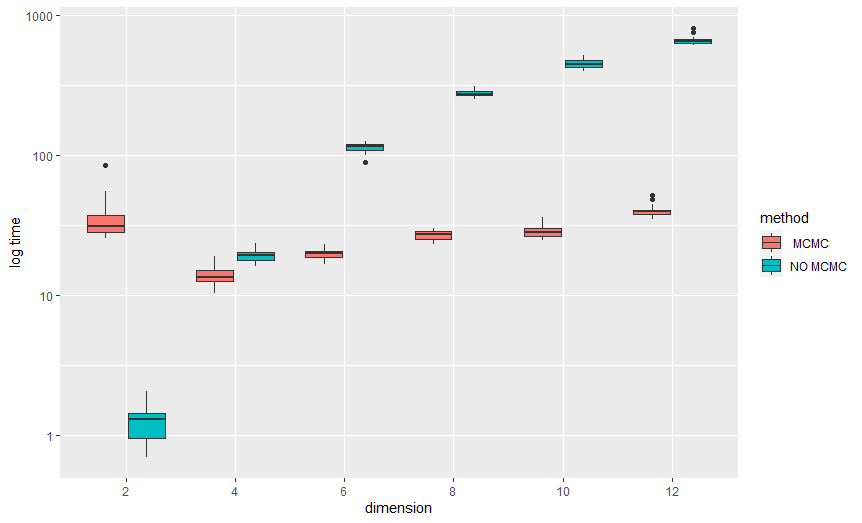}\\
\caption{Average calculation time (in seconds) as a function of size. "No MCMC" means "brute force Monte Carlo." Figure extracted from \cite{bergere2021}.}
\label{fig:timeMRM_courbes}
\end{figure}

We can then compare the theoretical value $p$ with bounds $(p^-_n,p^+_n)$ obtained through the nested algorithmic approach evoked above, after $n=200$ steps, and a basic one based on "brute force Monte Carlo". Relative precisions $(p^+_n-p^-_n)/p$ are provided on Table \ref{tab:table1.1} for several dimensions $d$, while Figure \ref{fig:timeMRM_courbes} displays how the computing time varies in function of $d$, for $p=5.10^{-4}$.

 \begin{table}[hbt!]
  \begin{center}
    \caption{Mean values of relative precision $(p^+_n-p^-_n)/p$ after 200 iterations for the MCMC-based  method.}
    \label{tab:table1.1}
        \vspace{0.5cm}
    
    \begin{tabular}{l c c  r} 
    \multicolumn{1}{c}{} &   \multicolumn{1}{c}{} & Brute force MC & MCMC\\\hline
    & \\
    \multicolumn{1}{c}{$d$} &   \multicolumn{1}{c}{$p$} &   & \\
    & \\
    \hline
    & \\
 \multirow{3}{*}{$2$ } & $5.10^{-2} $  & 0.233 & 0.208   \\       
                                 & $5.10^{-3} $   & 0.304 &  0.276  \\ 
                                 &$5.10^{-4} $    & 0.06 &  0.04 \\
                                 & \\ 
 \multirow{3}{*}{$3$ } & $5.10^{-2} $  & 1.18  & 1.09    \\       
                                 & $5.10^{-3} $   & 1.78 &  1.66  \\ 
                                 &$5.10^{-4} $    & 2.59 &  2.65  \\
                                 & \\ 
                                 
  \multirow{3}{*}{$4$ } & $5.10^{-2} $  & 2.53 &  2.38   \\       
                                 & $5.10^{-3} $   & 5.57 &  5.40  \\ 
                                 &$5.10^{-4} $    & 8.98 &  8.59  \\
                                 & \\ 
                                 \hline

    \end{tabular}
  \end{center}
\end{table}
  
Obviously, the computational time is significantly decreased using a MCMC approach, for the same order of magnitude for the highest dimensions. However, reaching small relative orders of magnitude remains difficult even for dimension 4, with $p^+_n \sim 5 p$, which fits with the results previously obtained by \cite{moutoussamy2013comparing}.  This remains a bit frustrating for the industrial practice, notably, but not hopeless. More simulations, and more clever simulations, for instance guided by a surrogate of $\Gamma$ sequentially updated, as the one proposed in \cite{bousquet2018approximation}, should be required to reach $p^+_n\sim p$.

\paragraph{Convexity and quasi-convexity.} It must be noticed that convexity and quasi-convexity properties on $g$ or $\Gamma$ can certainly help achieve this ambition. To our knowledge, the use of such properties has still little explored, apart in \cite{moutoussamy2015contributions}, and in the literature dedicated to computing probabilities related to PDE solving ; see for instance \cite{ullmann2024} for a summary presentation. Available results highlight speed orders on the distance between $\Gamma$ and surrogates of $\Gamma$, as exponential functions of the dimension; this distance can be bounded on convexity arguments \cite{ullmann2024}.  However,  bounding the implied relative error of the probability of failure still largely remains an open problem. \\

\section{Biasing surrogates to get conservative failure probabilities}\label{sec:trials}

\subsection{Principle}

When the use of a surrogate $\hat{g}_m$ trained from a $m-$sample is required for computational reasons, the estimation of $p$ by a $\hat{p}_{n,m}$ should be ideally such that
\begin{eqnarray}
p & \leq & \hat{p}_{n,m} \label{conservative.prob}
\end{eqnarray}
with a large probability $1-\alpha$, ideally with $\alpha=0$.  There are arguments to understand that $\alpha$ can possibly be strictly positive:  (a) the computations are made using a fixed data set (part of which is used to train the surrogate) ; it is therefore conceivable that not all the uncertainty about $\hat{g}_{m}(X)$ can be taken into account ; (b) when the dimension increases, as said in \cite{jakeman2022surrogate}, "surrogate modeling techniques are often built on lower dimensional subspaces identified by dimension-reduction techniques". These techniques introduce an error that is difficult to quantity. \\

To diminish $\alpha$ and more generally enforcing the so-called {\it first-order stochastic dominance constraint} \cite{jakeman2022surrogate}
\begin{eqnarray}
P\left(\hat{g}_{m}(X)< y\right) & \geq & P\left(g(X)< y\right) \ \ \ \forall y\in\R, \label{SDC}
\end{eqnarray}
an idea originated from \cite{Ducoffe2020} and also exploited by \cite{jakeman2022surrogate} is to "bias" (or "shift") the surrogate, replacing $x\mapsto\hat{g}_{m}(x)$ by
$$
x \mapsto\hat{g}_{m}(x) + \theta,
$$
to ensure conservativeness. Obviously, in our context $\theta$ should be chosen such that (\ref{conservative.prob}) is verified with $\hat{p}_{n,m}$ having the closest possible order of magnitude than the one of $p$. Obtaining too much difference in terms of order of magnitude should lead to modify the choice of the surrogate. 

\subsection{Using Bernstein concentration inequalities}

More precisely, \cite{Ducoffe2020} prove the following result, using a uniform Bernstein-type inequality:

\begin{theorem}[Inspired from \cite{Ducoffe2020}, Corollary 1]\label{theo.ducoffe}
Denote $\hat{g}_m$ a surrogate of $g$ trained over a $m-$sample and consider an independent test set of $n\geq 2$ iid values $X_1,\ldots,X_n$ drawn  over $\Omega$. Denote
\begin{eqnarray*}
\theta^* & = & \min\left(0,\min\limits_{1\leq i \leq n}\left\{\hat{g}_m(X_i)-g(X_i)\right\}\right).
\end{eqnarray*}
Then, with probability at least $1-\alpha$ over the choice of the test set, there exists a strictly positive constant $C<6$ such that 
\begin{eqnarray*}
P\left(\hat{g}_m(X)+\theta^* > g(X)\right) & \leq & B(n,\alpha):=\frac{C}{n}\log(n/\alpha).
\end{eqnarray*}
\end{theorem}

We can straightforwardly derive from it the following result, assuming $B(n,\alpha)=p$.

\begin{corollary}
Under the assumptions of Theorem \ref{theo.ducoffe},
\begin{eqnarray}
\Pp\left\{ P\left(\hat{g}_m(X)+\theta^* \leq g(X)\right) \geq 1-p\right\} & \geq & 1-\lambda(n,p) \label{control1}
\end{eqnarray}
with
\begin{eqnarray*}
\lambda(n,p) & = & \min\left(1,n\exp\left(-{np}/{C}\right)\right)
\end{eqnarray*}
Then, with probability at least $1-\lambda(n,p)$ over the choice of the test set,
\begin{eqnarray*}
\hat{p}_{m}:= P\left(\hat{g}_m(X)+\theta^* < y\right) & > & p.
\end{eqnarray*}
\end{corollary}

\begin{figure}[hbtp]
    \centering
    \includegraphics[scale=0.6]{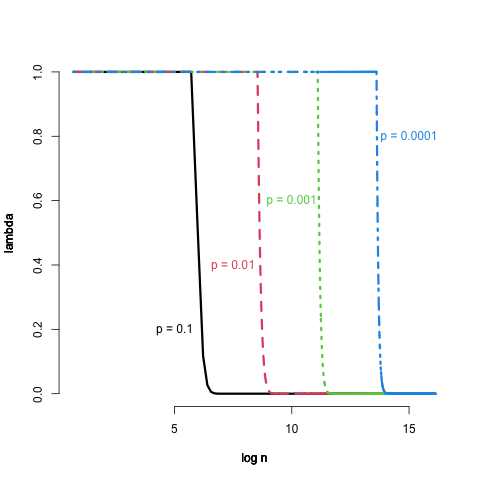}
    \caption{Values of $\lambda(n,p)$ with $C=6$, in function of $\log(n)$, for some values of $p$.}
    \label{fig:lambda}
\end{figure}

\begin{figure}[hbtp]
    \centering
    \includegraphics[scale=0.6]{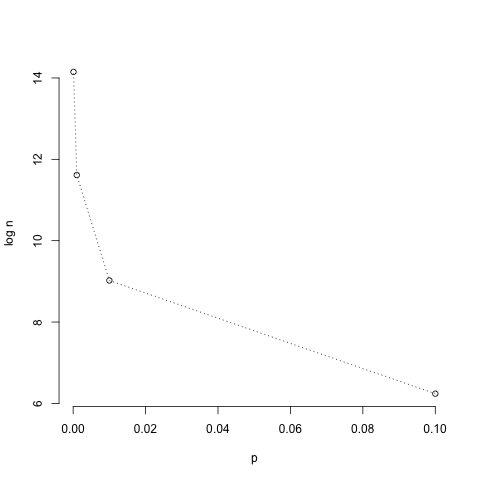}
    \caption{Values of $log(n)$ such that $\lambda(n,p)\simeq p$, for some typical values of $p$.}
    \label{fig:lambda2}
\end{figure}

While it is interesting to note that this result does not depend on the dimension of $X$, a large number of values $n$ is required to produce a nontrivial lower bound in (\ref{control1}), for  low probabilities $p$. This is exemplified by plotting representatives values of  $\lambda(n,p)$, setting $C=6$, in Figure \ref{fig:lambda}, and by computing the value of $n$ such that $\lambda(n,p)\sim p$, displayed (on the log scale) on Figure \ref{fig:lambda2}. Typically, one needs $n\simeq 512$ to obtain $\lambda(n,p)=p=10^{-1}$ then $n\simeq 8280$ to get  $\lambda(n,p)=p=10^{-2}$. In practice, this  cost remains probably too high for industrial applications. The authors \cite{Ducoffe2020} have logically proposed a complete procedure for learning $\theta$ in parallel of the surrogate parameters, which is enforced in the following numerical experiments.  


Reusing the toy example (\ref{toy.example.1}), numerical tests were conducted with neural networks surrogates. They were considered as good enough and retained for the study if both their accuracy and their  $Q^2$ predictivity coefficient \cite{fekhari2021model} were estimated above 90\% and 0.9, respectively. Neural networks were simple feedforward networks with logistic activation functions, 2 to 3 neuronal layers, each layers having 2 to 4 neurons. On Table \ref{tab:example2}, some results of these numerical tests are summarized. 

\begin{table}[hbtp]
\centering
  \caption{Some applied results on the reality of the conservative assessment of $p=P(g(X)<0)$ using simple neural networks, for the toy example (\ref{toy.example.1}), by the shifting/biasing approach originally proposed by \cite{Ducoffe2020} for creating conservative surrogates. The two last columns present the average conservative estimate $\hat{p}_n$ of $p$ produced by the surrogate, and the probability that a particular estimate actually be a wrong upper bound for $p$. These results are produced using $\lfloor100/(p(1-p))\rfloor$ repetitions of the procedure, and training datasets of length $50d$ to ensure a comparable precision of results.}
\label{tab:example2}
        \vspace{0.5cm}
    
\begin{tabular}{llllll}
\hline 
& & \multicolumn{4}{c}{Surrogate features}  \\
& \\
$p$ & $d$ & $Q^2$ without imposed bias & $Q^2$ with imposed bias & $\hat{p}_n$ & $\Pp(\hat{p}_n<p)$ \\
& \\
\hline 
& \\
$10^{-1}$ & 2 & 0.99 & 0.97 & $1.12.10^{-1}$ & 0.08 \\
          & 3 & 0.98 & 0.94 & $1.23.10^{-1}$ & 0.10 \\
          & 4 & 0.95 & 0.92 & $1.31.10^{-1}$ & 0.13 \\
          & 5 & 0.94 & 0.92 & $1.38.10^{-1}$ & 0.17 \\
          & \\
$10^{-2}$ & 2 & 0.98 & 0.97 & $1.9.10^{-2}$ & 0.09 \\
             & 3 & 0.96 & 0.94 & $2.4.10^{-2}$ & 0.14 \\
             & 4 & 0.94 & 0.91 & $3.1.10^{-2}$ & 0.21 \\
& \\
$10^{-3}$ & 2 & 0.96 & 0.94 & $2.8.10^{-3}$ & 0.18 \\
             & 3 & 0.95 & 0.91 & $3.9.10^{-3}$ & 0.25 \\
             & 4 & 0.95 & 0.90 & $6.2.10^{-3}$ & 0.34 \\

& \\
\hline 
\end{tabular}
\end{table}

These results confirm the theoretical reservations set out hereinbefore; it seems unacceptable to obtain, for "reasonable" low probabilities, upper bounds that can almost vary by an order of magnitude when the dimension remains low, and which turn out to be false upper bounds in a significant number of cases.

\subsection{Learning minimal bias}

A close approach but specifically  dedicated to obtaining upper bounds on probabilities was then proposed by \cite{jakeman2022surrogate}, inspired by the work \cite{viana2010using}, who build risk-averse surrogate models using stochastic dominance. Given a training sample $(x_i,y_i)_{1\leq i \leq m}$ and a parametric surrogate $x\mapsto g_{\eta}(x)$ mimicking $x\mapsto g(x)$,  they pose the problem of learning $(\theta,\eta)$ by choosing a specific weighted loss/cost function 
\begin{eqnarray}
\min\limits_{\theta,\eta} \sum\limits_{i=1}^m \omega_i\left(y_i - g_{\eta}(x_i)-\theta\right)^2 \label{loss.function}
\end{eqnarray}
to minimize under the first-order stochastic dominance constraints
\begin{eqnarray}
\left\{\begin{array}{ccc}\sum\limits_{i=1}^m \omega_i \1_{(-\infty,0]}\left(g_{\eta}(x_i)-g_{\eta}(x_1)\right) & \leq & \sum\limits_{i=1}^m \omega_i \1_{(-\infty,0]}\left(y_i - g_{\eta}(x_1)-\theta\right), \\
\ldots & \leq & \ldots \\
\sum\limits_{i=1}^m \omega_i \1_{(-\infty,0]}\left(g_{\eta}(x_i)-g_{\eta}(x_m)\right) & \leq & \sum\limits_{i=1}^m \omega_i \1_{(-\infty,0]}\left(y_i - g_{\eta}(x_m)-\theta\right)
\end{array}\right\}\label{hard.constraint}
\end{eqnarray}
This setting is the empirical version (ie., based on a given finite training dataset) of the estimation problem
\begin{eqnarray}
\left.
\begin{array}{l}
\min\limits_{\theta,\eta} \E\left[\left(Y-g_{\eta}(X)-\theta\right)^2\right] \\
\text{subject to $\theta+g_{\eta}(X) \succeq Y$}
\end{array}\right\}\label{theo.problem}
\end{eqnarray}
 where $Y_1 \succeq Y_2$ means that $Y_1$ dominates $Y_2$ with respect to the first stochastic order, ie. when
 \begin{eqnarray*}
P(Y_1\leq t) & \leq & P(Y_2\leq t)  \ \ \ \forall t\in\R.
 \end{eqnarray*}

Because of the discontinuity introduced by the indicator functions in (\ref{hard.constraint}), the authors\cite{jakeman2022surrogate} consider a continuous relaxation (originated from \cite{conti2018stochastic}) of the constrained problem (\ref{loss.function}-\ref{hard.constraint}), rewritten besides as a mixed integer optimization problem. Despite technical difficulties related to the choices of relaxations, the authors prove the relevance of this approximation and the conservativeness of the overall approach, ie.
\begin{eqnarray*}
P\left(\hat{g}_{m}(X) + \theta^* < y\right) & \geq & P\left(g(X)< y\right) \ \ \ \forall y\in\R,
\end{eqnarray*}
where $\hat{g}_{m}=g_{\eta^*}$, $(\eta^*,\theta^*)$ being the solution of (\ref{theo.problem}) and its empirical approximations. They also prove
its applicability for several uni- and multidimensional examples using polynomial chaos expansion, and show the feasibility of choosing $\theta$ dependent on $x$. In the same paper \cite{jakeman2022surrogate}, the authors defend another approach, based on the {\it risk quadrangle} \cite{rockafellar2013fundamental},  for building surrogates that conservatively estimate a specific risk measure associated with the subjective preferences of a stakeholder.

\section{Discussion}

This rapid review of tools and methods designed to support the construction and use of conservative estimators of probabilities of failure can be finally illustrated by a some typical numerical  results, displayed in Table \ref{summary.proba}.

\begin{table}[hbtp]
\centering 
\caption{Summary of some typical orders of magnitude of upper bounds for $p$ found in the dedicated literature ("un" means unpublished), considering toy models and more realistic examples. The value $n$ is the typical MC number of queries (training samples) $x\mapsto g(x)$ required for the computation (or samples for training a surrogate). The number of surrogate runs for a MC analysis is not precised, as it can be chosen as large as wished. Results $\Pp(p>p^+_n) $ are averaged using validation samples (from 100 to 9000 values) if surrogates are used. }
\label{summary.proba}

\vspace{0.5cm}

\begin{tabular}{ll}
\hspace{-1cm}
\begin{tabular}{cccccccc}
\hline
& \\
$p$ & $d$ & $n$ &  upper bound  & $\Pp(p>p^+_n) $ & use case & details & references \\
& & & $p^+_n$ & \\
& \\
\hline
& \\
\textcolor{blue}{$10^{-1}$} & 10 & 30 & \textcolor{blue}{$1.02.10^{-1}$} & 0 & wing weight & biased surrogate-based & \cite{jakeman2022surrogate} \\
& & & &&& (1D active space) & \\
\textcolor{blue}{$5.10^{-2}$} & 3 & 100 & \textcolor{blue}{$5.2.10^{-2}$} & 0 & toy model & biased surrogate-based & \cite{jakeman2022surrogate} \\
\textcolor{blue}{$5.10^{-2}$} & 10 & 100 & \textcolor{blue}{$10^{-1}$} & 0 & truss structure & biased surrogate-based & \cite{jakeman2022surrogate} \\ 
\textcolor{blue}{$2.1.10^{-3}$} & 1 & 32 & \textcolor{blue}{$2.5.10^{-3}$} & 0 & toy model & no surrogate,  & \cite{bernard2021recursive} \\
& & & & & & Lipschitz constant known & \\
\textcolor{blue}{$10^{-3}$} & 6 & 300 & \textcolor{blue}{$1.2.10^{-2}$} & 0 & toy model & no surrogate, $g$ monotonic & \cite{moutoussamy2013comparing} \\
\textcolor{blue}{$10^{-4}$} & 3 & 200 & \textcolor{blue}{$2.10^{-4}$} & 0 & toy model & no surrogate, $g$ monotonic & \cite{moutoussamy2013comparing} \\
\textcolor{blue}{$10^{-4}$} & 4 & 200 & \textcolor{blue}{$3.10^{-4}$} & 0 & hydraulics & no surrogate, $g$ monotonic & un. \\
\textcolor{blue}{$10^{-4}$} & 5 & 250 & \textcolor{blue}{$6.10^{-4}$} & 0 & toy model & no surrogate, $g$ monotonic & \cite{moutoussamy2013comparing} \\
& \\

\hline
\\
\end{tabular}
& 
\end{tabular}
\end{table}

These results, and the previous considerations, lead us to the following conclusions, which in turn allow us to outline a research program.

First, we only scratch the surface on how to bound failure probabilities with powerful methodologies involving the choice of surrogates. Indeed, targeted probabilities in research papers remains too "high" with respect to typical orders of magnitude tied to severe industrial risks. 

First, methodologies that are based on "biased" surrogates seem the most attractive to get non-asymptotic control of $p$, up to a given level of acceptability. While surrogates build on reduced bases usually come with bounds and could be ideal candidates, mixing constructive arguments proposed by \cite{jakeman2022surrogate} with other theoretical guarantees provided by concentration inequalities, as suggested by \cite{Ducoffe2020}, could be a relevant way to produce such guarantees. Refining the construction of such surrogates through sequential DOE, and obtaining concentration inequalities on the basis of martingale-type arguments in the spirit of De La Pe{\~n}a exponential inequalities \cite{bercu2015concentration}, seem a relevant avenue or research. Besides, ensuring the uniformity of first-order-stochastic dominance is probably not useful: we want to focus on some small subsets of $\Omega$.  

Second, it seems difficult, apart from some particular cases, to ensure strong geometrical assumptions directly on the behavior of $x\mapsto g(x)$. Perhaps the Lipschitz property seems more defensible than the monotonicity property in applications, but more probably, these two assumptions can be considered as true only for some part $X_{(1)}$ of the input $X$ conditionally to $X_{(2)}=X/X_{(1)}$. Provided the dimension of $X_{(1)}$ remains low, some small DOE could be used to provide conditional bounds on $p$ (given $X_{(2)}$). Other arguments to get bounds from the model $g$ itself could be inspired from the control techniques used on PDE solving by discretization (e.g., \cite{elfverson2016multilevel}. These techniques can provide hints to produce "cautious" surrogates, and offer an additional layer of guarantees for the end-user. Actually, the explicitly known or estimable geometric properties of these surrogates, such as their Lipschitz property and their monotonicity and convexity subdomains, make it possible to provide deterministic bounds on an estimator that is itself conservative and that would be produced from this surrogate, through the "biased surrogate" construction previously evoked.  

In this sense, an interesting research idea could be to produce some {additional criterion} to the usual ones used for assessing the relevance of surrogates (e.g., $Q^2$) to ensure a form of conservative bias in interesting subsets of $\Omega$. Such a work would fit with the recent development of a new criterion for Gaussian processes, interacting with the $Q^2$,  dedicated to ensure robust predictive properties of the surrogate \cite{marrel2023probabilistic}.  

In summary, the cases studied so far are usually very moderate in scale, and hardly reflect real cases (e.g. industrial, linked to critical systems), for which the need for a surrogate is particularly important. It therefore seems appropriate to launch a research program aimed at testing the scalability of the first methods based on first-order stochastic dominance, and to combine them with robust methods for calculating the bounds related to the geometry of these surrogates, through sequential explorations of the input space $\Omega$.

\section*{Acknowledgments}

The author thanks  Bastien Berg\`ere for having encoded and tested several ideas about monotonic Monte Carlo techniques, that have fed this review article.

\appendix
\section{Appendix: Semi-adaptive sampling and computing bounds within staircase subspaces}\label{appendix:MCMC}

As illustrated in Figure \ref{fig:mrm}, boundary algorithms based on monotonicity require to explore nested staircase subspaces (defined after each sampling batch $A\to A'$ by the space between dominated subspaces $(\UU^-(A'),\UU^+(A'))$). Basically, any strategy for producing a statistical estimator, computing and improving the bounds must be based on uniform sampling in this "tortured" nondominated space.  Rejection methods proposed in \cite{bousquet2012accelerated,bousquet2018approximation} unfortunately require an increasing number of simulations as the
as this space shrinks. A semi-adaptive Markov chain Monte Carlo method (MCMC) appears to be more efficient in dimension up to $d=15$. 

\begin{figure}[H] 
\centering
\includegraphics[scale=1.1]{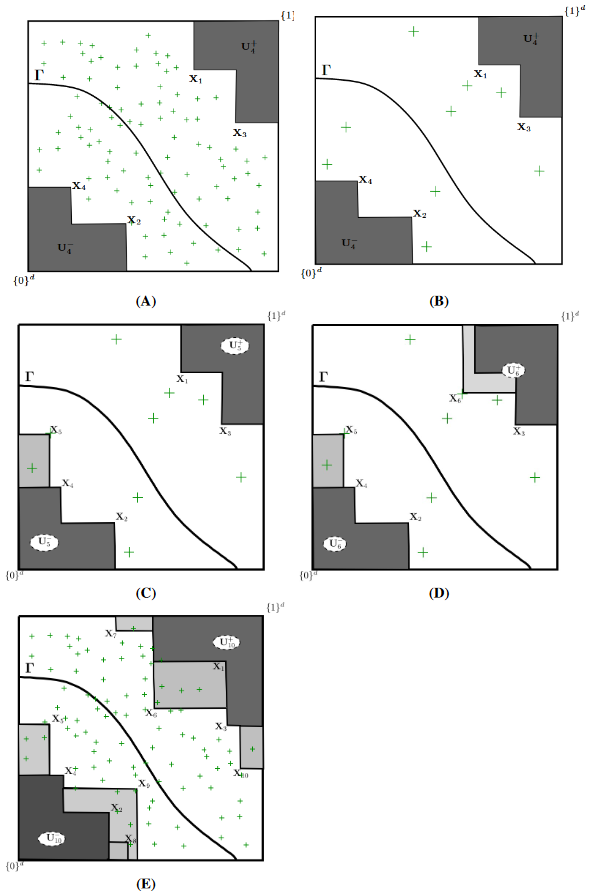}
\caption{Illustration in dimension 2 of the semi-adaptive MCMC method. From the previous chain which has converged (\textbf{A}), we obtain a quasi-iid sample following $\mathcal{U}(\UU_4)$ (\textbf{B}), which we use to generate the next realisations $(\bX_k)_{k=5,\ldots,10}$ following the uniform distribution on the respective spaces $(\UU_k)_{k=5,...,10}$ (\textbf{C}, \textbf{D}). The ergodic theorem is then applied to the entire trajectory of the chain to estimate the proportions of the sub-volumes $(\UU_k)_{k=5,...,10}$ (\textbf{E}).} 
\label{figMRMMCMC1} 
\end{figure}

\begin{figure}[H] 
\centering
\includegraphics[scale=1.1]{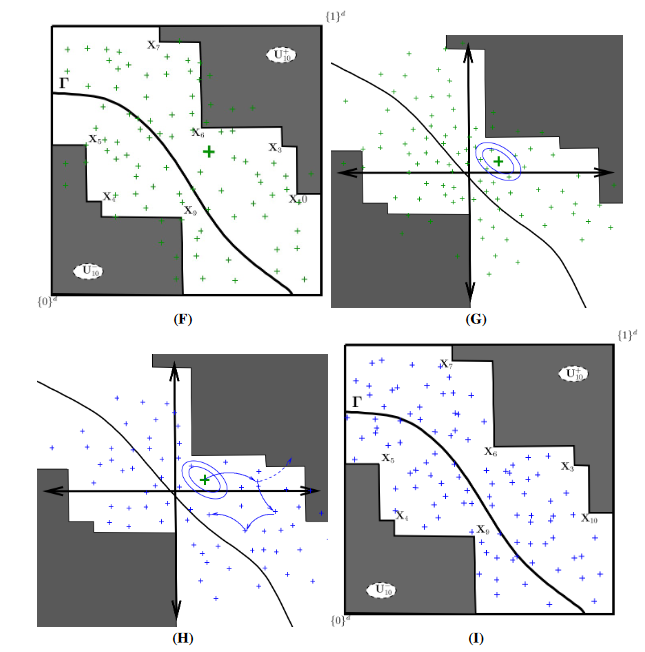}
\caption{Illustration in dimension 2 of the semi-adaptive MCMC method. After step  (\textbf{E}) on the previous figure, a point in the chain belonging to $\UU_{10}$ is chosen to initialise the next one (\textbf{F}). The Metropolis-Hastings random walk (\textbf{H}) is then evolved in the transformed space, taking Gaussian increments with a covariance equal to the empirical covariance of the previous chain (\textbf{G}). By inverse transformation we obtain a new chain converging to $\mathcal{U}(\UU_{10})$ (\textbf{I}).} 
\label{figMRMMCMC2} 
\end{figure}

More formally, assume that a $n-$DOE of queries $x\mapsto g(x)$ has been chosen in $\Omega=[0,1]^d$. This allows to define the two subspaces straightforwardly noted $\UU^-_n,\UU^+_n$, with $\UU_n=\Omega/(\UU^-_n\cup \UU^+_n)$. Our goal is to sample a $m-$sample from the uniform distribution with density 
\begin{eqnarray*}
f_n(x) & = & \frac{\1_{\{x\in \UU_n\}}}{p^+{n-1}-p^-{n-1}}.
\end{eqnarray*}
 Doing so allows to estimate the new bounds $(p^-_n,p^+_n)=(\lambda(\UU^-n),1-\lambda(\UU^+_n))$ and select new candidate vectors for  new queries $x\mapsto g(x)$. This can be done approximately using a MCMC technique, described beneath and illustrated on the successive graphs of Figures \ref{figMRMMCMC1}-\ref{figMRMMCMC2}, such that the computation of bounds is consistent thanks to the ergodic theorem. To explore efficiently the staircase subspace using a traditional random walk instrumental distribution, a transformation $\psi:\Omega\to\R^d$ of each possible design vector $x\in\UU_n$ is used, preserving the partial ordering $\preceq_P$: $\psi=(\phi^{-1},\ldots,\phi^{-1})$ where $\phi$ the cumulative distribution function of the standard reduced Gaussian distribution. \\
 
 More precisely, The semi-adaptiveness is defined in the following sense. Having constructed at a certain step $n$ a Metropolis-Hastings mechanism,  which has converged to the law $\mathcal{U}(\UU_n)$, its trajectory to estimate the volumes of interest and simulate uniformly at steps $n+1,\ldots, n+l$, for a fixed parameter $l<n$.A quasi-independent sample following $\mathcal{U}(\UU_n)$ is obtained by batch sampling. Then, from this sample, one simulates uniformly in the sub-regions of interest  by a simple accept-reject method.

\begin{algorithm}[H] \label{alg:MRM_MCMC}
 \DontPrintSemicolon
 \renewcommand{\algorithmcfname}{Algorithm}%

 \textbf{Initialization}\;
 Given a nontrivial staircase subspace $\UU_0$, sample $N$ values $\mathbf{X^{(0)}}\sim \mathcal{U}(\UU_0)$.  
 
 \Tq{n<N}{ 
 
 $\mathbf{1)}$ \textup{Compute $\mathbf{Z} = \psi(X_{(n-l)})$ and the empirical covariance matrix of transformed past $l-$trajectories \begin{equation*}
 \forall ~1\leq j,k \leq d~~~~ \hat{\Sigma}^{(n-l)}_{jk} =\frac{1}{N-1}\sum_{i=1}^{N}\left(  Z_{ij}-\bar{Z}_j \right)  \left( Z_{ik}-\bar{Z}_k \right)
 \end{equation*}
 then define } the innovation $\varepsilon  \sim \mathcal{N}_d(0, \alpha\hat{\Sigma}_{n-l} /d)$.  \\ 
 
 $\mathbf{2) \hspace{0.17cm} a)}$ Build the Markov chain $\mathbf{X}^{(n)}= (X^{(n)}_1,...,X^{(n)}_M)$ from a Metropolis-Hastings mechanism on the transformed variable $Z$, based on the stationary distribution $\mathcal{U}(\UU_n)$ and the transformed random walk proposal
 \begin{eqnarray*}
\tilde{X} & \overset{\cal{L}}{\sim} &  \psi^{-1}\left( \psi(X)+ \varepsilon\right).
 \end{eqnarray*}
 \hspace{0.5cm}$\mathbf{b)}$ Using usual decorrelation tools, this leads to a nearly independent batch sample  $\mathbf{\tilde{X}}^{(n)}$. \;
 \vspace{0.25cm}

 $\mathbf{3)}$ For each $k=1,,\ldots,l$, sample with rejection  $X_{n+k} \sim \mathcal{\UU}_{n+k}$ from $\mathbf{\tilde{X}}^{(n)}$.

 $\mathbf{4)}$ Estimate the volumes 
 \begin{eqnarray*}
 \widehat{\lambda(\UU_{n+k})} &= & \left( \frac{1}{M}  \sum_{i=m+1}^{m+M} \1_{\{\mathbf{Y}_i^{(n)} \in \UU_{n+k}\}} \right)  \widehat{\lambda(\UU_{n})} \\
     \widehat{\lambda(\UU^+_{n+k})} &=& \left(\frac{1}{M}  \sum_{i=m+1}^{m+M} \1_{\{\mathbf{X}_i^{(n)} \in \hat{\UU}^+_{n+k}\}}\right) \widehat{\lambda(\UU_{n})}
 \end{eqnarray*}
 where, $ \forall k \in \llbracket 1,\ldots, l\rrbracket, ~~ \forall i:$
 \begin{eqnarray*}
\1_{\{\mathbf{Y}_i^{(n)} \in \UU_{n+k}\}} & = & \left[ \1_{\{\mathbf{Y}_i^{(n)} \succeq X_{n+k}\}}\1_{\{\mathbf{X}_{n+k}\in \UU-} + \1_{\mathbf{Y}_i^{(n)} \preceq X_{n+k}\}}\1_{\{\mathbf{X}_{n+k}\in \UU+\}}\right] \1_{\{\mathbf{Y}_i^{(n)} \in \UU_{n+k-1}\}}.  
 \end{eqnarray*}\;
 $\mathbf{5)}$ Update the bounds and update $n \longleftarrow n +l$
 }
 \caption{Semi-adaptive MCMC for estimating probability bounds.}
\end{algorithm}

\bibliographystyle{plain}
\bibliography{biblio}

\end{document}